\documentclass{siamltex}

\usepackage{amsfonts} \usepackage{amsmath} \usepackage{amssymb} \usepackage{epsfig}

\newcommand{\bR}{\mathbb{R}} \newcommand{\bC}{\mathbb{C}}
\newcommand{\bZ}{\mathbb{Z}} \newcommand{\real}{{\rm Re}\,} \newcommand{\imag}{{\rm
Im}\,}  \newcommand{\arccosh}{\rm arccosh}
\newcommand{\eps}{\epsilon} \newcommand{\veps}{\varepsilon}

\newtheorem{teorema}{Theorem} 
\newtheorem{ilustracion}{Illustration} \newtheorem{nota}{Remark}

\def \precision { \varepsilon }

\title{A spectral order method for inverting sectorial Laplace transforms}

\author{Mar\'ia L\'opez-Fern\'andez \footnotemark[2] \and Cesar Palencia
\footnotemark[2] \and Achim Sch\"adle \footnotemark[5]}

\begin{document}

\maketitle

\renewcommand{\thefootnote}{\fnsymbol{footnote}}

\footnotetext[2]{Departamento de Matem\'atica Aplicada, Universidad de Valladolid,
Valladolid, Spain.  ~E-mail: {\tt \{marial, palencia\}@mac.cie.uva.es}. Supported by
DGI-MCYT under project MTM 2004-07194 cofinanced by FEDER funds.}

\footnotetext[5]{ZIB Berlin, Takustr.~7, D-14195 Berlin, Germany.  ~E-mail: {\tt
schaedle@zib.de}. Supported by the DFG Research Center \textsc{Matheon} "Mathematics
for key technologies" in Berlin.}

\begin{abstract}
Laplace transforms which admit a holomorphic extension to some sector strictly
containing the right half plane and exhibiting a potential behavior are considered. A
spectral order, parallelizable method for their numerical inversion is proposed.  The
method takes into account the available information about the errors arising in the
evaluations. Several numerical illustrations are provided.
\end{abstract}


\section{Introduction}\label{introduccion}

In a variety of situations, the problem arises of inverting numerically the Laplace
transform $U(z)$ of a given mapping of interest $u(t)$. Roughly speaking, it turns
out that the wider the set $W$ where $U(z)$ can be computed is, the easier the
inversion results. For instance, if $W$ is an interval $(a,b)$ then the numerical
inversion becomes an ill-posed problem \cite{Ang,Cun,Cunh}. On the other hand, if $W$
is the complement of some bounded region, then the efficient Talbot's method
\cite{Riz,Tal} is at hand.

In the present paper we focus on the particular situation where $W$ is a sector
symmetric with respect to the real axis, strictly containing the right half plane,
and we assume that $U(z)$ exhibits a potential behavior on $W$. We say then that
$U(z)$ is {\it sectorial}. Precisely, there is a renewed interest in the numerical
inversion of {\it sectorial} mappings \cite{Gav,Gavr,LoP04,Mcl,She}, mainly due to
its applicability to linear, non-homogeneous evolution equations of parabolic type
(both in the context of abstract IVP's and Volterra equations), as well as their
discretizations in space \cite{Ash,Bak}. Notice that the applicability of the
inversion approach, in the {\it sectorial} setting, demands in practice that the
source term of the parabolic equation must be approximated efficiently (at least
locally) by holomorphic mappings \cite{Gavr,Mcl}. This difficulty is overcome in
\cite{LoLuPaSch,SchLoLu}, where the ideas in the present paper are adapted so as to
provide accurate reconstructions of the traditional Runge-Kutta approximations to the
solutions of such parabolic problems. These reconstructions require no regularity on
the source term of the problem.

In the present paper we consider the issue of the numerical inversion of {\it
sectorial} mappings by itself. To fix ideas, let
$$ u:(0, +\infty) \to X
$$ be a locally integrable mapping, taking values in a Banach space $X$, with
exponential growth. Denote by
$$ U(z)=\int_0^{+\infty}e^{-zt} u(t)\, dt
$$ its Laplace transform. We will always assume that $U(z)$ admits a holomorphic
extension to the complement $W$ of some acute sector
\begin{equation}\label{sector}
  \Sigma_{\delta}=\{z\in\bC : |\arg (-z)|\leq \delta\}, \qquad
  0<\delta<\frac{\pi}{2},
\end{equation}
an that there exist constants $M > 0$ and $ \mu\in \bR$ such that
\begin{equation}\label{desU}
  \|U(z)\|\leq\frac{M}{|z|^{\mu}}, \qquad z \notin \Sigma_\delta.
\end{equation}
The last requirement, with $\mu \ge 1$, means that $u$ admits a bounded and
holomorphic extension to any sector of the form $|\arg(z)| \le \delta'$, with $0 <
\delta' < \pi/2-\delta$. If $\mu < 1$, we select an integer number $m \ge 1$, with
$m+\mu \ge 1$, and set $V(z) = U(z)/z^m$. Then, by the previous remark, $V(z)$ is the
Laplace transform of a mapping $v :(0,+\infty) \to X$, which admits a bounded and
holomorphic extension to sectors with semi-angle $\delta'$ as before, and now $u$ is
understood to be the derivative of order $m$ of $v$.

Notice that in case $U(z)$ satisfies a similar inequality
$$ \|U(z)\|\leq\frac{M}{|z-\omega|^{\mu}}, \qquad z\notin \omega + \Sigma_{\delta},
$$ for some $\omega \in \bR$, then, by using the shifting theorem, the inversion of
$U(z)$ is reduced to the one of a Laplace transform $\widetilde U(z)$ fulfilling
(\ref{desU}). Since the respective originals $u(t)$ and $\widetilde u(t)$ are related
by $u(t) = e^{\omega t} \widetilde u(t)$, then we can just approximate $\widetilde
u(t)$. This is why the analysis is restricted to the situation $\omega = 0$, i.e. to
(\ref{desU}).

The goal is to numerically reconstruct $u$ from knowledge of a moderate number of
evaluations of $U(z)$ at suitable nodes $z \notin \Sigma_\delta$. Let us point out
that, from a practical point of view, it is essential to take into account that these
evaluations are going to be affected by errors.

The starting point of the method we propose is the well-known inversion formula
\begin{equation}\label{invLT}
  u(t)=\frac{1}{2\pi i}\int_{\Gamma} e^{tz}U(z)\, dz,\qquad t>0,
\end{equation}
where $\Gamma$ is a suitable path connecting $-i\infty$ to $+i\infty$ which, in our
setting, can be taken so as to guarantee the absolute convergence of the integral
appearing in (\ref{invLT}). As in \cite{Gav,LoP04,Mcl}, we choose $\Gamma$ the branch
of a hyperbola pa\-ra\-me\-tri\-zed by a mapping $S:(-\infty, +\infty) \to \bC $
admitting a holomorphic extension to a horizontal strip around the real axis. The
numerical method we propose is simply the truncated trapezoidal rule, applied to the
definite integral arising after parametrizing (\ref{invLT}) by $S$, used with $2n+1$
nodes $x_k = kh$, $-n \le k \le n$, and a suitable step size $h>0$. The properties of
$S$ allow us to use the ideas and results in \cite{Ste,Sten}, where the trapezoidal
rule applied to holomorphic mappings on strips is considered. Let us comment that the
fast decay of our integrand \cite{Gav,LoP04} yields an improvement of the more
general estimates in \cite{Ste,Sten}.

Very often, for instance in the context of IVP's (see Illustration~\ref{ilustracion3}
in Section~\ref{secexperimentos}), the main computational effort of the method is due
to the evaluations of $U(z)$ at the nodes $z_k = S(x_k)$, $-n \le k \le n$. An
important feature of the present approach is that the same evaluations can be used to
approximate $u(t)$ at different $t > 0$ \cite{LoP04,Riz}. Accordingly, our goal is to
obtain a uniform error estimate for the approximation of $u(t)$ on intervals of the
form $[t_0, \Lambda t_0]$, with given $t_0 >0$ and $\Lambda \ge 1$, rather than at a
fixed $t >0$. Essentially, this was the aim in \cite{LoP04}, whose basic estimates we
borrow. Notice also that the algorithm presents two levels of parallelism since,
first, the evaluations of $U(z)$ at the involved nodes and, second, the evaluations
of $u(t)$ at a selected finite set of values of $t \in [t_0,\Lambda t_0]$, can be
carried out on different processors.

In the present paper, by considering a different choice of the geometrical and scale
parameters from the one in \cite{LoP04}, we improve the results there in two
different ways: \renewcommand{\labelenumi}{(\roman{enumi})}
\begin{enumerate}
\item We get a better error bound, which now turns out to be a genuine spectral
  estimate of the form $O(e^{-cn})$.
\item We also get a weaker dependence of the exponential factor $c$ on $\Lambda$,
  since now $c = O(1/ \ln \Lambda)$.
\end{enumerate}
This means, in practice, that with a moderate number of evaluations of $U(z)$ we can
accurately approximate $u(t)$, uniformly on intervals $[t_0,t_1]$ with $\Lambda =
t_1/t_0 >> 1$, let us say $\Lambda = 50$.

On the other hand, for the choice of parameters we propose, the precision $\rho$ used
in the evaluations of $U(z)$ at the required nodes plays a more relevant role than in
\cite{LoP04}. In fact, ignoring that we always have $\rho > 0$ would result in large
actual errors for $n >> 1$, as simple numerical experiments show (see Illustration~1
in Section~\ref{secexperimentos}). This drawback is overcome by minimizing the
estimate we get for the actual error (Theorem~\ref{teorema2}), which leads to a
$(\rho,n)$-dependent choice of parameters. With this choice, the actual error finally
behaves for moderate $n$ like $O(e^{-cn})$, with $c = O(1/\ln \Lambda$), and for
large $n$ like $O(\rho)$.  This optimal choice of parameters demands, of course, some
information about the size of $\rho$. In the absence of it, we propose an
$n$-dependent choice of parameters for which the actual error behaves like
$O(\rho+e^{-cn})$, with $c = O(1/(\ln n + \ln \Lambda))$. All the above estimates are
uniform on $t_0 \le t \le \Lambda t_0$, with fixed $t_0>0$ and $\Lambda
>1$. Moreover, the error constants are made explicit in the analysis and turn out to
be reasonable.

The outline of the paper is as follows. In Section~\ref{secmetodonumerico} we
describe the numerical method and show, in Theorem~\ref{teorema1}, how to achieve (i)
and (ii).  The propagation of errors is studied in Section~\ref{secactualerror}. The
choice of parameters is considered in Section~\ref{secparametros} and four simple
numerical illustrations of the theoretical results are provided in
Section~\ref{secexperimentos}.


\section{The numerical method}\label{secmetodonumerico}
Given $\delta$ in (\ref{sector}) and following the ideas in \cite{LoP04}, we select
$\alpha,\, d > 0$ such that
\begin{equation}\label{condicionangulos}
  0<\alpha -d < \alpha + d < \frac\pi 2 -\delta~.
\end{equation}
Defining
\begin{equation}
  \label{eq:Tpre}
  T(w) = 1-\sin(\alpha+iw)
\end{equation}
this mapping transforms each horizontal straight line $ \imag w = y,\ -d \leq y \leq
d, $ into the left branch of the hyperbola given by
\begin{equation}
  \label{eq:hyperbola}
  \bigg(\frac{\real z - 1}{\sin(\alpha-y)} \bigg)^2 - \bigg(\frac{\imag z}{\cos
  (\alpha-y)}\bigg)^2 = 1,
\end{equation}
with center at $(1,0)$, foci at $(0,0)$ and $(2,0)$, whose asymptotes make angles
$\pm [\pi/2-(\alpha-y)]$ with the real axis.  Therefore, $T$ transforms the
horizontal strip
$$ D_d=\{ z\in \bC : |\imag z|\leq d \}
$$ into the region in the complex plane limited by the left branches corresponding to
$y = \pm d$ in \eqref{eq:hyperbola}.

Introducing a parameter $\lambda > 0$, the parametrization of $\Gamma$ in
(\ref{invLT}) can be defined as
$$ \Gamma =\{ \lambda T(x) : x\in \bR \},
$$ i.e. $\Gamma$ is the branch of a hyperbola corresponding to the image of the real
axis under $S=\lambda T$. This results in
$$ u(t)=\int_{-\infty}^{+\infty}G_t(x)\,dx~, \qquad t>0,
$$ where $G_t: D_d \to X,\ t>0$, is the mapping
$$ G_t(w)=-\frac{\lambda}{2\pi i}\exp(\lambda tT(w))U(\lambda T(w))T'(w).
$$ Once the parameters $\alpha$, $d$, and $\lambda$ have been fixed, we set $x_k=kh,\
k\in \bZ$, and consider the approximation to $u(t)$ given by
\begin{equation}\label{invnumerica}
  u_n(t)=h\sum_{k=-n}^{n}G_t(x_k),\qquad t>0.
\end{equation}
The proof of the main result in \cite{LoP04} (Theorem~2), shows that for $\mu=1$ in
(\ref{desU})
\begin{equation}\label{cotaantigua}
\|u(t)-u_n(t)\| \le M \cdot \varphi(\alpha,d) \cdot L(\lambda t\sin(\alpha-d)) \cdot
e^{\lambda t}\Big(\frac{1}{e^{2\pi d/h}-1}+\frac{1}{e^{\lambda t \sin\alpha
\cosh(nh)}}\Big),
\end{equation}
where
$$ \varphi(\alpha,d) = \frac{2}{\pi}\sqrt{\frac{1+\sin(\alpha+d)}{1-\sin(\alpha+d)}},
$$ and $L(x)$, $x > 0$, is the function
$$ L(x)=1+|\ln(1-e^{-x})|.
$$ Notice that $L(x)$ is decreasing in $x$, $L(x)\approx |\ln x|$, for $x \to 0^{+}$
and $L(x)$ tends to $1$, for $x\to +\infty$.

As we commented in the Introduction, in many applications the computational effort to
obtain $u_n(t)$ is mainly due to the evaluations of $U(z)$ at $z=\lambda T(x_k),\ -n
\leq k \leq n$, but these evaluations could be carried out in parallel. Another
attractive feature of (\ref{invnumerica}) is that the same evaluations of $U(z)$ can
be used to compute $u_n(t)$ for different $t > 0$. In fact, as we see below, with the
appropriate choice of parameters, we can use the same evaluations of $U(z)$ so as to
have a spectral estimate
$$ \|u(t)-u_n(t)\| = \mathcal{O}(e^{-cn}),
$$ uniform on intervals $t_0 \leq t \leq t_1$. The exponential factor $c$ turns out
to depend weakly on the ratio $\Lambda = t_1/t_0$, given that $ c = \mathcal{O} ( 1/
\ln \Lambda)$.

For simplicity the next theorem is restricted to the situation $\mu = 1$ in
(\ref{desU}). The cases $\mu > 1$ and $\mu < 1$ are treated in subsequent remarks.
\begin{teorema}\label{teorema1}
  Assume that $U$ satisfies (\ref{desU}) with $\mu=1$.  Fixing $\alpha$ and $d$
  according to (\ref{condicionangulos}), for $t_0 > 0$, $ \Lambda \geq 1$, $0 <
  \theta < 1$ and $n\geq 1$, the following choice of parameters
  \begin{equation}\label{parametros}
    h=\frac{1}{n}a(\theta),\qquad \lambda = \frac{2\pi d n (1-\theta)}{t_0 \Lambda
      a(\theta)},
  \end{equation}
  with
  $$ a(\theta)=\arccosh \Big( \frac{\Lambda}{(1-\theta)\sin \alpha} \Big),
  $$ yields the uniform estimate on $t_0 \leq t \leq \Lambda t_0$
  \begin{equation}\label{cotanueva}
    \|u(t)-u_n(t)\| \leq M \cdot \varphi (\alpha,d) \cdot L( \lambda t_0
    \sin(\alpha-d))\cdot \frac{2\epsilon_n(\theta)^{\theta}}{1-\eps_n(\theta)},
  \end{equation}
  where
  $$ \eps_n(\theta)=\exp\Big( -\frac{2\pi d }{a(\theta)}n\Big).
  $$
\end{teorema}
The theorem shows, just by selecting any $0 < \theta < 1$, a genuine spectral order
of convergence in $n$ of the form $O(e^{-cn})$, where $c=O(1/\ln \Lambda)$
(cf. \cite{Gav,LoP04}).

{\em Proof}. Set $\sigma =\lambda t_0$. For $t_0 \leq t \leq \Lambda t_0$,
(\ref{cotaantigua}) implies the uniform bound
$$ \|u(t)-u_n(t)\| \leq M \cdot \varphi (\alpha,d)\cdot L(\sigma \sin(\alpha
-d))\cdot e^{\Lambda \sigma} \Big( \frac{1}{e^{2\pi d/h}-1} +\frac{1}{e^{\sigma \sin
\alpha \cosh (nh)}} \Big).
$$ Our choice of $h$ and $\lambda$ is precisely the one guaranteeing that
$$ \exp\Big( \frac{2\pi d}{h} \Big) = \exp ( \sigma \sin \alpha \cosh(nh)) =
\frac{1}{\eps_n(\theta)},
$$ hence
$$ \frac{1}{e^{2\pi d/h}-1} +\frac{1}{e^{\sigma \sin \alpha \cosh (nh)}} \leq \frac
{2e^{-2\pi d/h}}{1-e^{-2\pi d/h}} = \frac {2\eps_n(\theta)}{1-\eps_n(\theta)}.
$$ The proof ends after remarking that
$$ e^{\Lambda \sigma}\eps_n(\theta) = \eps_n(\theta)^{\theta-1}\eps_n(\theta) =
\eps_n(\theta)^{\theta}.\qquad \endproof
$$

To end the section we comment, in the two following remarks, on the situation $\mu
\ne 1$ in (\ref{desU}). We omit details in the proofs, which are completely analogous
to the one of Theorem~\ref{teorema1}.

\begin{nota}
  \label{nota1}
  {\rm Assume that $U$ satisfies (\ref{desU}) with $\mu > 1$.  By Remark~1 in
    \cite{LoP04} we have
    $$ \|u(t)-u_n(t)\| \le M \cdot \varphi(\alpha,d,\mu) \cdot L(\lambda
    t\sin(\alpha-d)) \cdot \frac{e^{\lambda t}}{\lambda^{\mu-1}}\Big(\frac{1}{e^{2\pi
    d/h}-1}+\frac{1}{e^{\lambda t \sin\alpha \cosh(nh)}}\Big),
    $$ where
    $$\varphi(\alpha,d,\mu) =
    \frac{2}{\pi}\sqrt{\frac{1+\sin(\alpha+d)}{(1-\sin(\alpha+d))^{2\mu-1}}}.
    $$ Thus, for $0 < \theta < 1$, the same choice of values for $h$ and $\lambda$ as
    in Theorem~\ref{teorema1} gives the bound
    $$ \|u(t)-u_n(t)\| \leq M \cdot \varphi (\alpha,d,\mu) \cdot L( \lambda t_0
    \sin(\alpha-d)) \cdot \lambda^{1-\mu} \cdot
    \frac{2\epsilon_n(\theta)^{\theta}}{1-\eps_n(\theta)},
    $$ uniformly for $t_0 \le t \le \Lambda t_0$. This estimate is again spectral in
    $n$, since
    $$ \lambda^{1-\mu} = O\Big( \Big( \frac{\Lambda t_0}{n} \Big)^{\mu-1} \Big).
$$}
\end{nota}

\begin{nota}
  \label{nota2}
  {\rm Assume now that $U$ satisfies (\ref{desU}) with $\mu<1$.  By Remark~1 in
    \cite{LoP04}, for a fixed $s \in (0,1)$, there holds
    $$ \|u(t)-u_n(t)\| \le M \cdot \varphi_s(\alpha,d,\mu) \cdot L(s\lambda t
    \sin(\alpha-d)) \cdot \frac{e^{\lambda t}}{t^{1-\mu}}\Big(\frac{1}{e^{2\pi
    d/h}-1}+\frac{1}{e^{s\lambda t \sin\alpha \cosh(nh)}}\Big),
    $$ where now
    $$\varphi_s(\alpha,d,\mu) =
    \frac{2}{\pi}\sqrt{\frac{1+\sin(\alpha+d)}{1-\sin(\alpha+d)}}
    \Big(\frac{1-\mu}{(1-s) e \sin(\alpha-d)} \Big)^{1-\mu}.
    $$ In this situation, for $\theta \in (0,1)$ we choose
    $$ h=\frac{1}{n}a_s (\theta),\qquad \lambda = \frac{2\pi d n (1-\theta)}{t_0
    \Lambda a_s(\theta)},
    $$ where
    $$ a_s(\theta )=\arccosh \Big (\frac{\Lambda}{s(1-\theta)\sin \alpha }\Big).
    $$ Setting
    $$ \eps_{s,n}(\theta)=\exp \Big( \frac{-2 \pi d n }{a_s(\theta)}\Big),
    $$ we get the spectral estimate
    $$ \|u(t)-u_n(t)\| \leq M \cdot \varphi_s (\alpha,d,\mu)\cdot L( s\lambda t_0
    \sin(\alpha-d))\cdot t_0^{\mu-1}
    \frac{2\eps_{s,n}(\theta)^{\theta}}{1-\eps_{s,n}(\theta)},
    $$ uniformly for $t_0 \le t \le \Lambda t_0$. }

\end{nota}

\section{Error propagation}\label{secactualerror}
Numerical experiments (see Section \ref{secexperimentos}), show that for large values
of $n$ the estimate~(\ref{cotanueva}) is not longer true in practice. The explanation
of this apparently contradictory behavior lays in the influence of the errors when
evaluating $U$ and the elementary functions involved. For the sake of simplicity, we
consider first the case $\mu=1$ in (\ref{desU}).  The situations $\mu > 1$ and $\mu <
1$ are considered in subsequent remarks.

Let $z_k = \lambda T(x_k),\ -n \leq k \leq n$, be the nodes used in
(\ref{invnumerica}). Clearly, in practice, as numerical approximation to $u(t)$ we
actually obtain
\begin{equation}\label{invnumericareal}
  \bar{u}_n(t) = \sum_{k=-n}^{n} \omega_k(t) U_k ,
\end{equation}
where, for $-n \le k \le n$, $\omega_k(t) \in \bC$ and $U_k \in X$, are
approximations to
$$ -\frac{\lambda h}{ 2 \pi i} \exp(\lambda t z_k )T'(x_k),
$$ and $U(z_k)$, respectively.

To estimate the actual error $\|u(t)-\bar{u}_n(t)\|$ we need to make some assumptions
on the approximations used. To this end, we are going to focus on two frequent
possibilities, depending on whether we have information on absolute or relative
errors due to the evaluations. To be precise, we are going to assume that there
exists $\rho >0$ such that, simultaneously for all $-n \leq k \leq n$, we have either
\begin{equation}\label{abserrorevaluacion}
  \| U(z_k) - U_k \| \le \rho \quad \mbox{and} \quad \omega_k(t) =-\frac{\lambda h}{
  2 \pi i} \exp(\lambda t z_k )T'(x_k)
\end{equation}
or
\begin{equation}\label{relerrorevaluacion}
  \|\exp(\lambda t z_k)T'(x_k)U(z_k) - \omega_k(t) U_k \| \leq \rho \|\exp(\lambda t
  z_k)T'(x_k)U(z_k)\|.
\end{equation}

Situation (\ref{abserrorevaluacion}) arises for instance when $U_k \approx U(z_k) $
are provided by means of some auxiliary routine, let us say by solving a linear
system, with prescribed accuracy $\rho$ and moreover the errors due to the
evaluations of the elementary functions involved turn out to be negligible compared
to $\rho$. Situation (\ref{relerrorevaluacion}) is typical when $U(z)$ is an
elementary function.

The next theorem yields an estimate of the actual error for these situations.  We
maintain the notation introduced in Theorem~\ref{teorema1}.

\begin{teorema}\label{teorema2}
  Assume that $U$ satisfies (\ref{desU}) with $\mu=1$. Fix $\alpha$, $d$ according to
  (\ref{condicionangulos}). For $t_0 >0$, $ \Lambda \geq 1,\ 0 < \theta < 1$ and $n
  \geq 1$, select the parameters
  $$ h=\frac{1}{n}a(\theta),\qquad \lambda = \frac{2\pi d n (1-\theta)}{t_0 \Lambda
  a(\theta)}.
  $$ Assume also that $\omega_k(t) \in \bC,\ t_0 \le t \le t_1$, $U_k \in X$, $-n\leq
  k \leq n$, satisfy either (\ref{abserrorevaluacion}) or (\ref{relerrorevaluacion}).
  Then, the actual error is estimated by
  \begin{equation}\label{cotareal}
    \|u(t)- \bar{u}_n(t)\| \leq M \cdot \Phi \cdot Q \cdot \bigg( \precision
    \eps_n(\theta)^{\theta-1} +\frac{\eps_n(\theta)^{\theta}}{1-\eps_n(\theta)}
    \bigg),
  \end{equation}
  uniformly on $t_0 \leq t \leq \Lambda t_0$, where either
  \renewcommand{\labelenumi}{(\alph{enumi})}
  \begin{enumerate}
  \item $\precision=\rho /(M t_0)$,
    $$ \Phi=\max
    \bigg\{\frac{2}{\pi}\sqrt{\frac{1+\sin(\alpha+d)}{1-\sin(\alpha+d)}},
    \frac{1}{\pi e \sin \alpha} \bigg\}
    $$ and
    $$ Q = \max \bigg \{ 2L(\lambda t_0 \sin(\alpha-d)), \frac{\ln n}{\ln n-1}
    \bigg[\frac{\ln n}{2n}+L\Big(\frac{\lambda t_0 \sin\alpha}{\ln n}\Big) \bigg]
    \bigg \}
    $$ in case (\ref{abserrorevaluacion}) holds, or
  \item $\precision = \rho$,
    $$ \Phi =\frac{2}{\pi}\sqrt{\frac{1+\sin(\alpha+d)}{1-\sin(\alpha+d)}}
    $$ and
    $$ Q = \max \{ 2L(\lambda t_0 \sin(\alpha-d)), 1/2 (h + L(\lambda t_0 \sin
    \alpha))\},
    $$ in case (\ref{relerrorevaluacion}) holds.
  \end{enumerate}
\end{teorema}

Notice that $Q$ depends logarithmically on $\alpha,\, d, \, 1-\theta$ and $\Lambda$.

The estimate (\ref{cotareal}) given by the theorem, with a fixed $0 < \theta < 1$,
shows again a spectral order of convergence $O(e^{-cn})$, with $c=O(1/\ln \Lambda)$,
but only for moderate n, to be more precise, as long as $\eps_n(\theta) \geq
\veps$. In fact, for fixed $\theta$, (\ref{cotareal}) goes to $+\infty$ as $n \to
+\infty$. However, this apparent drawback is overcome by selecting $\theta$ in a
suitable way, as we explain in Section 4.

\begin{proof}
By writing
$$ \|u(t)-\bar{u}_n(t)\| \leq \|u(t)-u_n(t)\| + \|u_n(t)-\bar{u}_n(t)\|,
$$ and noticing that, for the corresponding $Q$, (\ref{cotanueva}) implies
$$ \|u(t)- u_n(t)\| \leq M\cdot\Phi\cdot Q
\frac{\eps_n(\theta)^{\theta}}{1-\eps_n(\theta)},
$$ the proof is reduced to show that
\begin{equation}
  \label{reduc}
  \|u_n(t)-\bar{u}_n(t)\| \le M\cdot\Phi\cdot Q \precision \eps_n(\theta)^{\theta-1}.
\end{equation}

Assume first that (\ref{abserrorevaluacion}) holds. This situation was already
studied in \cite{LoP04}, where it was proved that
\begin{equation}\label{propagacionerrorabs}
  \|u_n(t)-\bar{u}_n(t)\| \leq \frac{\rho \ln n}{2\pi e (\ln n-1) \sin\alpha}
    \frac{e^{\lambda t}}{t} \bigg[\frac{\ln n}{n}+2L\Big(\frac{\lambda t
    \sin\alpha}{\ln n}\Big) \bigg],
\end{equation}
whence, after recalling that $\precision = \rho/(t_0 M)$ and noticing that
\begin{equation}
  \label{noteteo2}
  e^{\Lambda \lambda t_0} = \eps_{n} (\theta)^{\theta-1},
\end{equation}
we readily obtain (\ref{reduc}).

Assume now that (\ref{relerrorevaluacion}) holds. Proceeding as in the proof of
Lemma~1 and Theorem~2 in \cite{LoP04}, and denoting
$$ \varphi (\alpha,0) = \frac{2}{\pi}\sqrt{\frac{1+\sin\alpha}{1-\sin\alpha}},
$$ we get
\begin{eqnarray*}
  \|u_n(t)-\bar{u}_n(t)\| &\leq& \frac{\rho M e^{\lambda t}}{2 \pi} h \sum_{k=-n}^{n}
  e^{-\lambda t \sin \alpha \cosh x_k} \bigg| \frac{T'(x_k)}{T(x_k)} \bigg| \\ [5pt]
  \nonumber &\leq& \frac{M \varphi (\alpha,0)}{4} \rho e^{\lambda t}h \sum_{k=-n}^{n}
  e^{-\lambda t \sin \alpha \cosh x_k} \\ [5pt] \nonumber &\leq& \frac{M \varphi
  (\alpha,0)}{2}\rho e^{\lambda t} \bigg( h + \int_0^{+\infty} e^{-\lambda t \sin
  \alpha \cosh x}\, dx \bigg) \\ [5pt] \nonumber &\leq& \frac{M \varphi
  (\alpha,0)}{2} \rho e^{\lambda t}( h + L(\lambda t \sin \alpha) ).
\end{eqnarray*}
Hence, using again (\ref{noteteo2}) and the inequality $\varphi (\alpha,0) \le
\varphi (\alpha,d)$, we deduce (\ref{reduc}).  \qquad \end{proof}

The behavior $\mu \ne 1$ in (\ref{desU}) is considered in the following remarks,
whose proofs are a combination of Remark~1, Remark~2 and the arguments used in the
proof of Theorem~\ref{teorema2} in \cite{LoP04}. Notice that
(\ref{propagacionerrorabs}) is independent of $\mu$.

\begin{nota}
  \label{nota3}
  {\rm Assume that $\mu >1$ in (\ref{desU}) and fix $0 < \theta <1$. Then, for the
    choice of parameters in Theorem~\ref{teorema2} and uniformly on $t_0 \le t \le
    \Lambda t_0$, we have: \renewcommand{\labelenumi}{(\alph{enumi})}
    \begin{enumerate}
    \item in case (\ref{abserrorevaluacion}) holds
      $$ \|u(t)- \bar{u}_n(t)\| \leq M \cdot \Phi \cdot Q \cdot \bigg( \precision
      \eps_n(\theta)^{\theta-1} + \lambda^{1-\mu}
      \frac{\eps_n(\theta)^{\theta}}{1-\eps_n(\theta)} \bigg)
      $$ with $\precision=\rho/(Mt_0)$,
      $$ \Phi=\max
      \bigg\{\frac{2}{\pi}\sqrt{\frac{1+\sin(\alpha+d)}{(1-\sin(\alpha+d))^{2\mu-1}}},
      \frac{1}{\pi e \sin \alpha} \bigg\}
      $$ and
      $$ Q = \max \bigg \{ 2L(\lambda t_0 \sin(\alpha-d)), \frac{\ln n}{\ln n-1}
      \bigg[\frac{\ln n}{2n}+L\Big(\frac{\lambda t_0 \sin\alpha}{\ln n}\Big) \bigg]
      \bigg \},
      $$
    \item in case (\ref{relerrorevaluacion}) holds
      $$ \|u(t)- \bar{u}_n(t)\| \leq M\cdot \Phi\cdot Q\cdot \lambda^{1-\mu}
      \cdot\bigg( \precision \eps_n(\theta)^{\theta-1}
      +\frac{\eps_n(\theta)^{\theta}}{1-\eps_n(\theta)} \bigg),
      $$ with $\precision=\rho$,
      $$ \Phi
      =\frac{2}{\pi}\sqrt{\frac{1+\sin(\alpha+d)}{(1-\sin(\alpha+d))^{2\mu-1}}}
      $$ and
      $$ Q = \max \{ 2L(\lambda t_0 \sin(\alpha-d)), 1/2( h + L(\lambda t_0 \sin
      \alpha))\}.
      $$
    \end{enumerate}
  }
\end{nota}

\begin{nota}\label{nota4}
  {\rm Assume that $\mu < 1$ in (\ref{desU}) and fix $0 < s, \theta <1$. Then, for
    the choice of parameters in Remark~\ref{nota2} and uniformly on $t_0 \le t \le
    \Lambda t_0$, we have:

    \renewcommand{\labelenumi}{(\alph{enumi})}
    \begin{enumerate}
    \item in case (\ref{abserrorevaluacion}) holds
      $$ \|u(t)- \bar{u}_n(t)\| \leq M\cdot \Phi\cdot Q\cdot \bigg( \precision
      \eps_{s,n}(\theta)^{\theta-1} + t_0^{\mu-1}
      \frac{\eps_{s,n}(\theta)^{\theta}}{1-\eps_{s,n}(\theta)} \bigg),
      $$ with $\precision=\rho/(Mt_0)$,
      $$ \Phi=\max
      \bigg\{\frac{2}{\pi}\sqrt{\frac{1+\sin(\alpha+d)}{1-\sin(\alpha+d)}}\Big(
      \frac{1-\mu}{(1-s) e\sin(\alpha-d)}\Big)^{1-\mu}, \frac{1}{\pi e \sin \alpha}
      \bigg\}
      $$ and
      $$ Q = \max \bigg \{ 2L(s\lambda t_0 \sin(\alpha-d)), \frac{\ln n}{\ln n-1}
      \bigg[\frac{\ln n}{2n}+L\Big(\frac{\lambda t_0 \sin\alpha}{\ln n}\Big) \bigg]
      \bigg \},
      $$
    \item in case (\ref{relerrorevaluacion}) holds
      $$ \|u(t)- \bar{u}_n(t)\| \leq M\cdot \Phi\cdot Q\cdot \bigg( \lambda^{1-\mu}
      \precision \eps_{s,n}(\theta)^{\theta-1} + t_0^{\mu-1}
      \frac{\eps_{s,n}(\theta)^{\theta}}{1-\eps_{s,n}(\theta)} \bigg),
      $$ with $\precision=\rho$,
      $$ \Phi =\frac{2}{\pi}\sqrt{\frac{1+\sin(\alpha+d)}{1-\sin(\alpha+d)}} \Big(
      \frac{1-\mu}{(1-s) e\sin(\alpha-d)}\Big)^{1-\mu},
      $$ and
      $$ Q = \max \{ 2L(s\lambda t_0 \sin(\alpha-d)), 1/2(h + L(s \lambda t_0 \sin
      \alpha))\}.
      $$
    \end{enumerate}
  }

\end{nota}

\section{The choice of parameters}\label{secparametros}

With Theorem~\ref{teorema2} in mind, we now try to derive a strategy for the choice
of parameters. First of all, (\ref{cotareal}) shows that it is of interest to select
$\alpha$ away from zero and $\alpha+d$ away from $\pi/2$.  The dependence of the
actual error on $\alpha-d$ is less important, since it is logarithmic.

Suppose $\alpha$ and $d$ have been already chosen, then for a given $n$ we take $h$
and $\lambda$ as indicated in Theorem~\ref{teorema2} and we fix $0 < \theta <
1$. Assume also that we have an estimation of $\rho$ and set $\varepsilon=\rho
/(Mt_0)$ or $\veps=\rho$ as in Theorem 2. Then, since in practice we always have
$\rho > 0$ and hence $\varepsilon > 0$, it turns out that $\precision
\eps_n(\theta)^{\theta-1} \to +\infty$, as $n \to +\infty$. Hence, it is clear that
increasing the number of nodes might result in a worse estimate (\ref{cotareal}). In
fact, increasing $n$ may result in worse approximations, as
Illustration~\ref{ilustracion1} in Section~\ref{secexperimentos} shows.

To overcome this drawback we let $\theta$ be a free parameter for the moment. Given
$\precision >0$ and $n$, after selecting $\alpha$ and $d$, neglecting the logarithmic
factor $Q$ and taking into account that typically $\eps_n(\theta) << 1$, the best
thing we can do is to choose $0 < \theta <1$ so as to minimize the term
\begin{equation}\label{fminimizar}
  \precision \eps_n(\theta)^{\theta-1} + \eps_n(\theta)^\theta,
\end{equation}
i.e. we must tune $\theta$ depending on $\precision >0$ and $n$.  By a direct
calculation it can be proven that the first derivative of $\eps_n(\theta)^{\theta-1}$
with respect to $\theta$ is increasing in $\theta$. The same is true for
$\eps_n(\theta)^\theta$ (in this case the proof, though elementary, is more
difficult). We conclude that the expression in (\ref{fminimizar}) is a convex
function of $\theta$. Moreover, its limit either for $\theta \to 0+$ or $\theta \to
1-$ is $+\infty$.  Therefore, (\ref{fminimizar}) attains its minimum exactly for one
value $\theta_{\varepsilon,n} \in (0,1)$, which is the one we propose to be
used. Though it is not easy to express the dependence of $\theta_{\veps,n}$ on $n$
and $\precision$, this can be easily done numerically (see Section 5).

Since, up to logarithmic factors, the choice $\theta=\theta_{\veps,n}$ in
(\ref{cotareal}) is optimal, it is clear that with this choice we get for the actual
error: \renewcommand{\labelenumi}{(\alph{enumi})}
\begin{enumerate}
\item A spectral order of convergence $O(e^{-cn})$ with $c=O(1/\ln \Lambda)$, for
moderate values of $n$, since this is true for any value $0< \theta <1$.

\item The errors are not propagated. In fact, already with the non-optimal choice
$$ \theta=1-\frac 1 n,
$$ (\ref{cotareal}) reads
\begin{equation}\label{cotanooptimal}
\|u(t)-\bar{u}_n(t)\|=O(\veps+e^{-cn}),
\end{equation}
uniformly on $t_0 \le t \le \Lambda t_0$, with $c=O(1/(\ln \Lambda +\ln n))$. This
remark tells us that, for large values of $n$, the actual error saturates at level
$\veps$, as observed in the numerical experiments (see Section 5).
\end{enumerate}

In the previous discussion it was essential to assume that we had some information
about $\veps$. Notice that, even in case we do not have such an information, the
choice $\theta =1-1/n$, which led to (\ref{cotanooptimal}), is always available. This
bound is almost spectral in $n$, depends weakly on $\Lambda$ and prevents error
amplification.

\section{Numerical illustrations}\label{secexperimentos}
In this section we give four numerical illustrations. The first two ones concern
elementary Laplace transforms which are assumed to be computed with a relative error
of order $\rho \approx \mbox{eps}$, where eps stands for the machine precision
($\mbox{eps}=10^{-16}$ in our computations). In the last two illustrations we do not
assume any information about the errors due to the computations of the Laplace
transforms.

\begin{ilustracion}\label{ilustracion1}
  {\rm We first show by means of a simple example, that for $n>>1$ (\ref{cotanueva})
    fails in the presence of errors in the evaluations. To this end, we consider the
    mapping $u(t)=e^{-t}$, whose Laplace transform is $U(z)=1/(1+z)$.

    This function satisfies (\ref{desU}) for all $\delta > 0$ and $M = 1/\sin
    \delta$. We fix $\theta =0.5,\ \alpha=0.7,\ d=0.6$ and choose the parameters $h,\
    \lambda$ as stated in the theorem for all the values of $n$. In
    Fig.~\ref{figerrprop} we plot in a semilogarithmic scale the absolute actual
    error, i.e.
    $$\ln \max_{t\in [t_{0}, \Lambda t_{0}]}\|u(t)-\bar{u}_n(t)\|$$ versus $n$
    (recall that $\bar{u}_n(t)$ stands for the actual computed approximation to
    $u(t)$, see (\ref{invnumericareal})). This is done for $\Lambda=5,50$ and
    $t_0=1$. This figure shows that the error decays exponentially for the first
    values of $n$, saturates near $\veps$ level and then grows like $O(e^{cn})$.
    \begin{figure}
      \begin{center}
        \epsfig{file=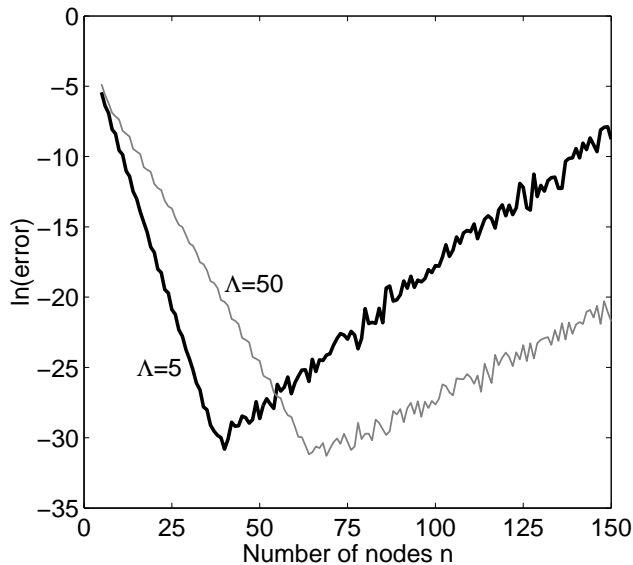, height=3in}
      \end{center}
      \caption{$\ln \max_{t\in [t_{0}, \Lambda t_{0}]}\|u(t)-\bar{u}_n(t)\|$
        versus $n$ for $u$ in Illustration~\ref{ilustracion1}, with $\theta =0.5$
        fixed, $\alpha=0.7,\ d=0.6$ and $t_0=1$. The gray line corresponds to
        $\Lambda = 50$ and the black one to $\Lambda=5$.} \label{figerrprop}
    \end{figure}

    Next we tune parameters as explained in Section 4.  For $\Lambda =5,\ 50$, in
    Fig.~\ref{figoneoverzplusone} (left) we plot the optimal values of $\theta$
    against $n$. In Fig.~\ref{figoneoverzplusone} (right) we plot
    $$ \ln \max_{t \in [t_0, \Lambda t_0]} \|u(t)-\bar{u}_n(t)\|
    $$ (continuous line) and the logarithm of the corresponding values of the
    theoretical error estimate (dashed line) obtained in Theorem~\ref{teorema2},
    versus $n$, once $\theta$ is optimal.  We maintain $\alpha=0.7,\ d=0.6$ and
    $t_0=1$.

    \begin{figure}
      \begin{center}
        \epsfig{file=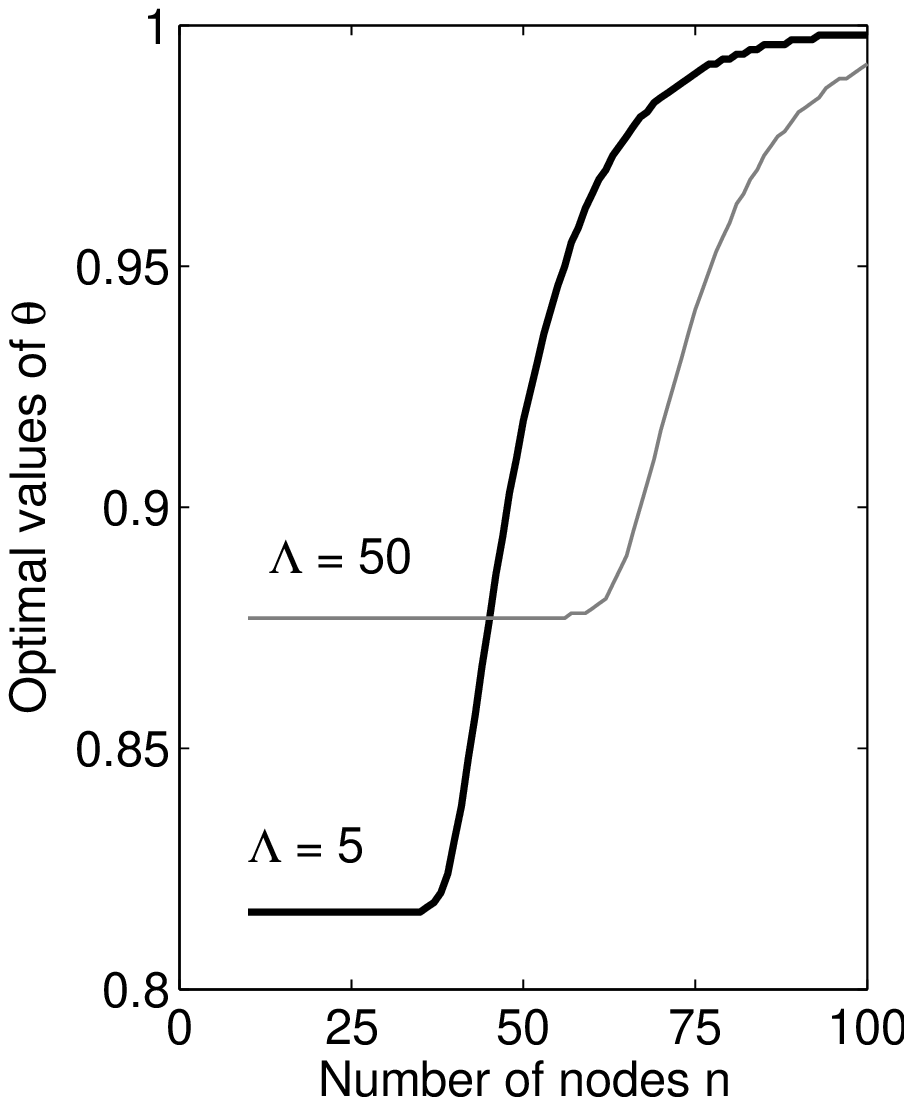, height=3in} \epsfig{file=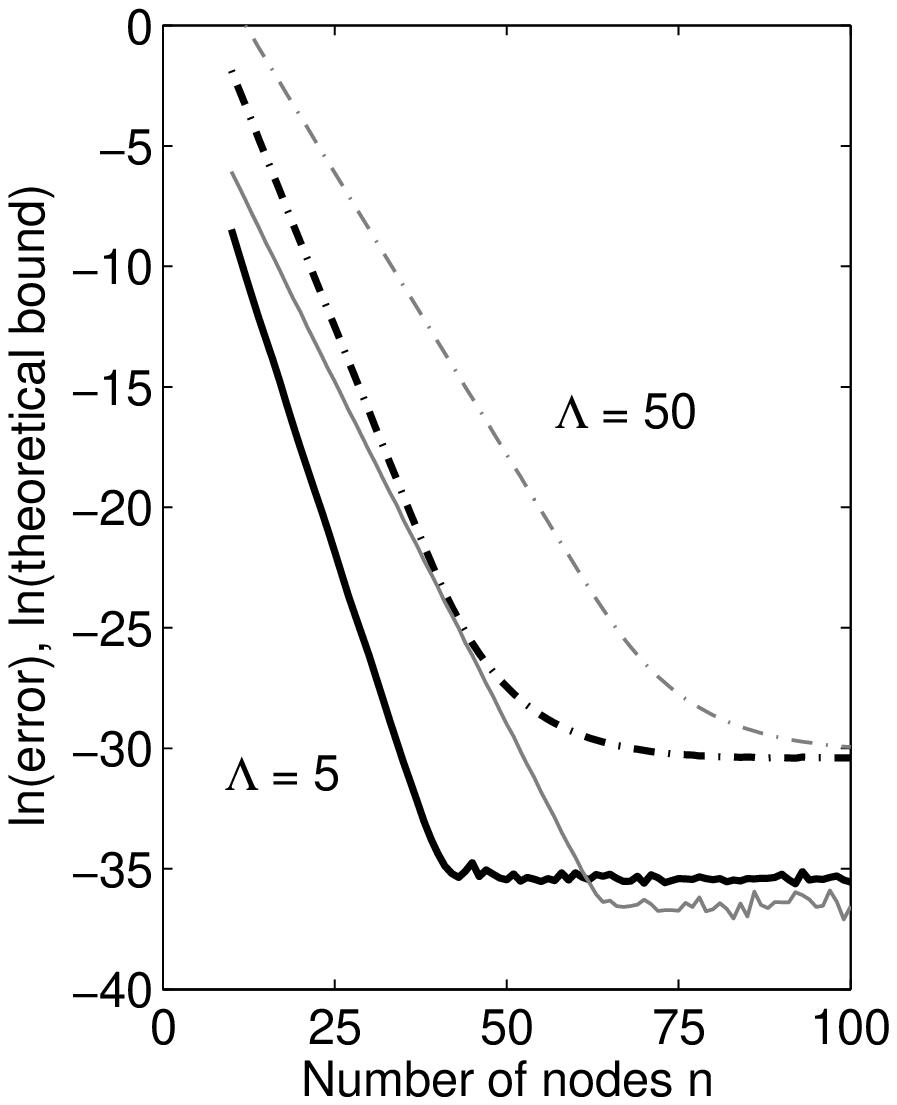,
        height=3in}
      \end{center}
      \caption{Left: Optimal $\theta$ versus $n$. Right: Natural logarithms
      of $\max_{t \in [t_0, \Lambda t_0]}\|u(t)- \bar{u}_n(t)\|$ (continuous) and the
      theoretical estimate (dashed) versus $n$, for $u$ in
      Illustration~\ref{ilustracion1}. The gray lines correspond to $\Lambda = 50$
      and the black ones to $\Lambda=5$.}
      \label{figoneoverzplusone}
    \end{figure}
  }
\end{ilustracion}

\begin{ilustracion}\label{ilustracion2}
  {\rm Take $\beta=1.5$ and set
    $$ U(z)=\frac{z^{\beta-1}}{z^{\beta}+1},
    $$ i.e., $U(z)$ is the Laplace transform of
$$ u(t)=M_{\beta}(-t^{\beta}),
$$ where $M_\beta$ stands for the Mittag-Leffler function of order $\beta$ (see
\cite{Pol}).  Notice that $U$ satisfies (\ref{desU}) for any $\delta \in
(\pi/3,\pi/2)$, with $\mu=1$ and $M=1/\sin(\beta(\pi-\delta))$. We consider here as
exact solution the one computed with 500 nodes and take $\alpha=\pi/12,\ d=0.25$ and
$t_0=1$.

    This example was already considered in \cite{LoP04}. In order to compare the
    performance of the strategy proposed in \cite{LoP04} with the one proposed in the
    present paper, we first compute $\bar{u}_n(t)$ by selecting the parameters as in
    \cite{LoP04}. In Fig.~\ref{figolderror} (left) we plot in semilogarithmic scale
    the theoretical estimate and actual errors for $\Lambda =2,5$, which are
    acceptable. In Fig.~\ref{figolderror} (right) we do the same for $\Lambda=50$ and
    conclude that the approach in \cite{LoP04} is not at all useful for large values
    of $\Lambda$. However, the corresponding computation by using the strategy in
    Section~4, yields the plot in Fig.~\ref{figMLeffler}, which shows a satisfactory
    spectral order of convergence even for $\Lambda=50$.

    \begin{figure}[ht]
      \begin{center}
        \epsfig{file=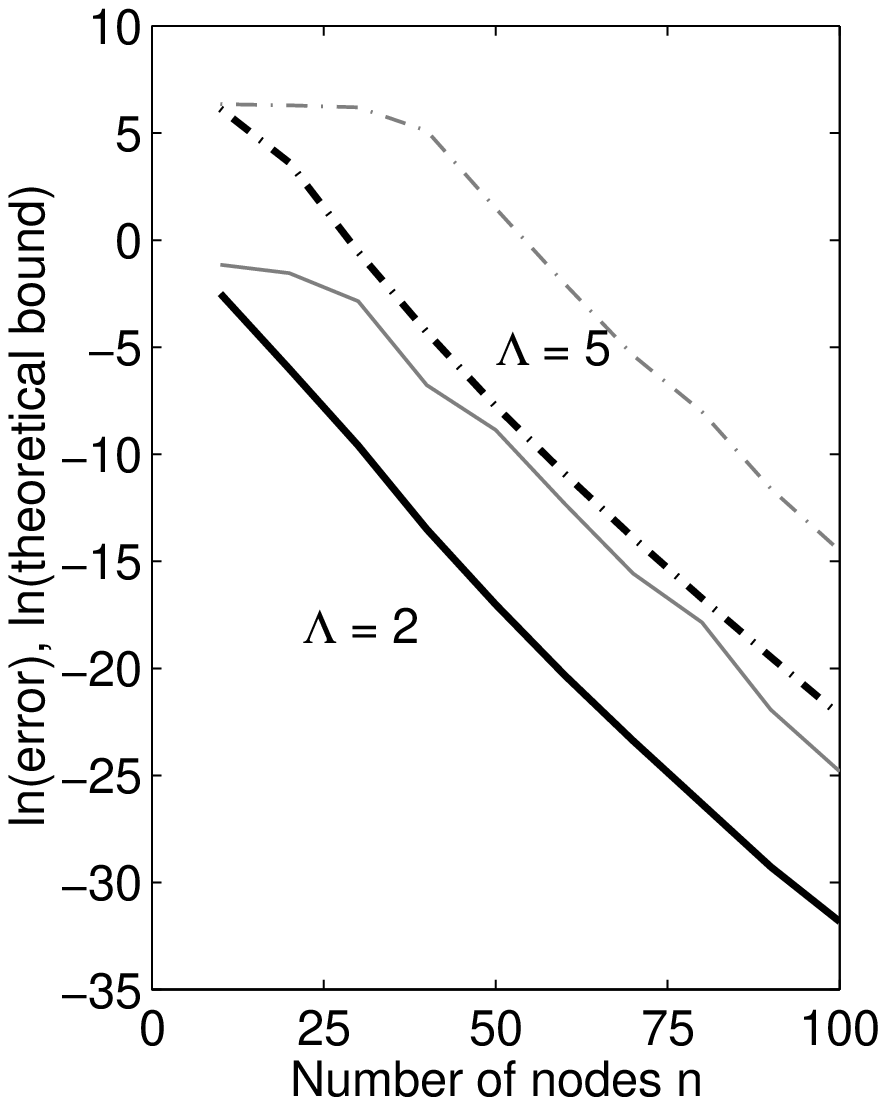, height=3in} \epsfig{file=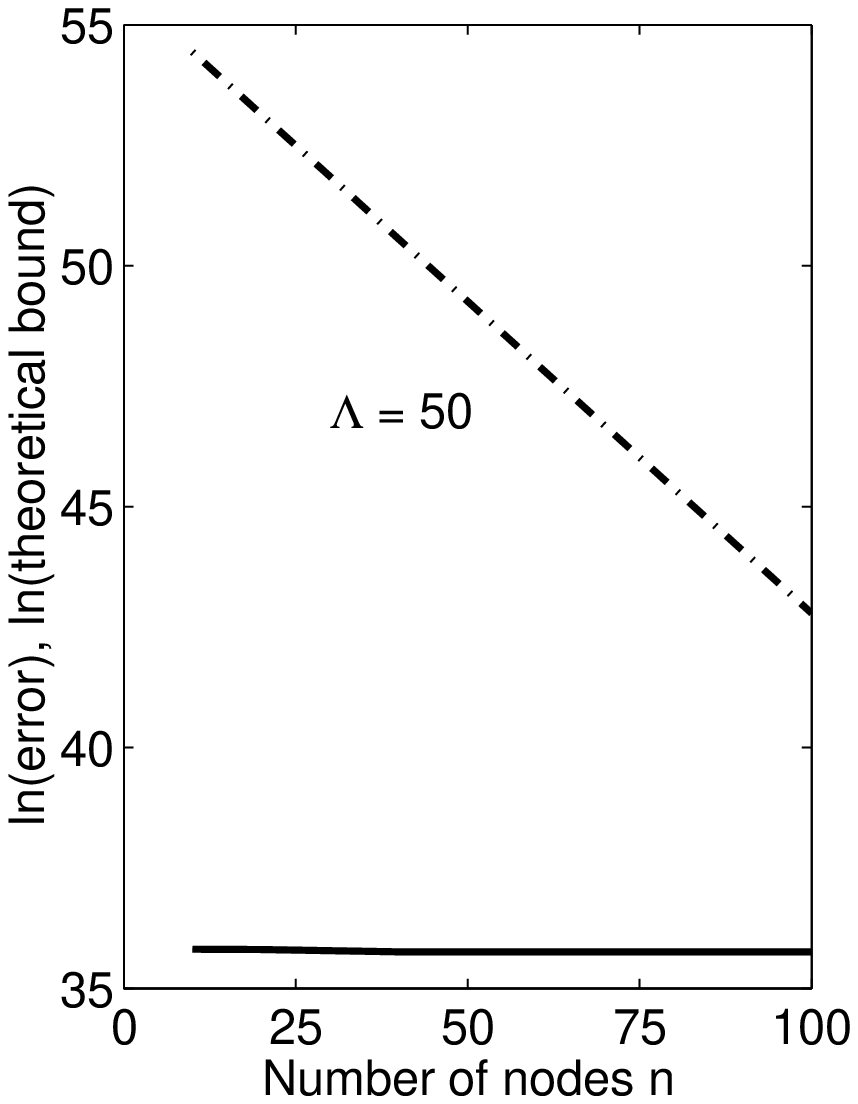, height=3in}
      \end{center}
      \caption{Natural logarithms of $\max_{t \in [t_0, \Lambda t_0]}\|u(t)-\bar{u}_n(t)\|$
      (continuous) and the theoretical estimate (dashed) versus $n$, for $u$ in
      Illustration~\ref{ilustracion2} proceeding as in \cite{LoP04} for
      $\delta=\pi/3$, $t_0=1$.  The gray lines correspond to $\Lambda = 50$ and the
      black ones to $\Lambda=5$.}
      \label{figolderror}
    \end{figure}

    \begin{figure}[ht]
      \begin{center}
        \epsfig{file=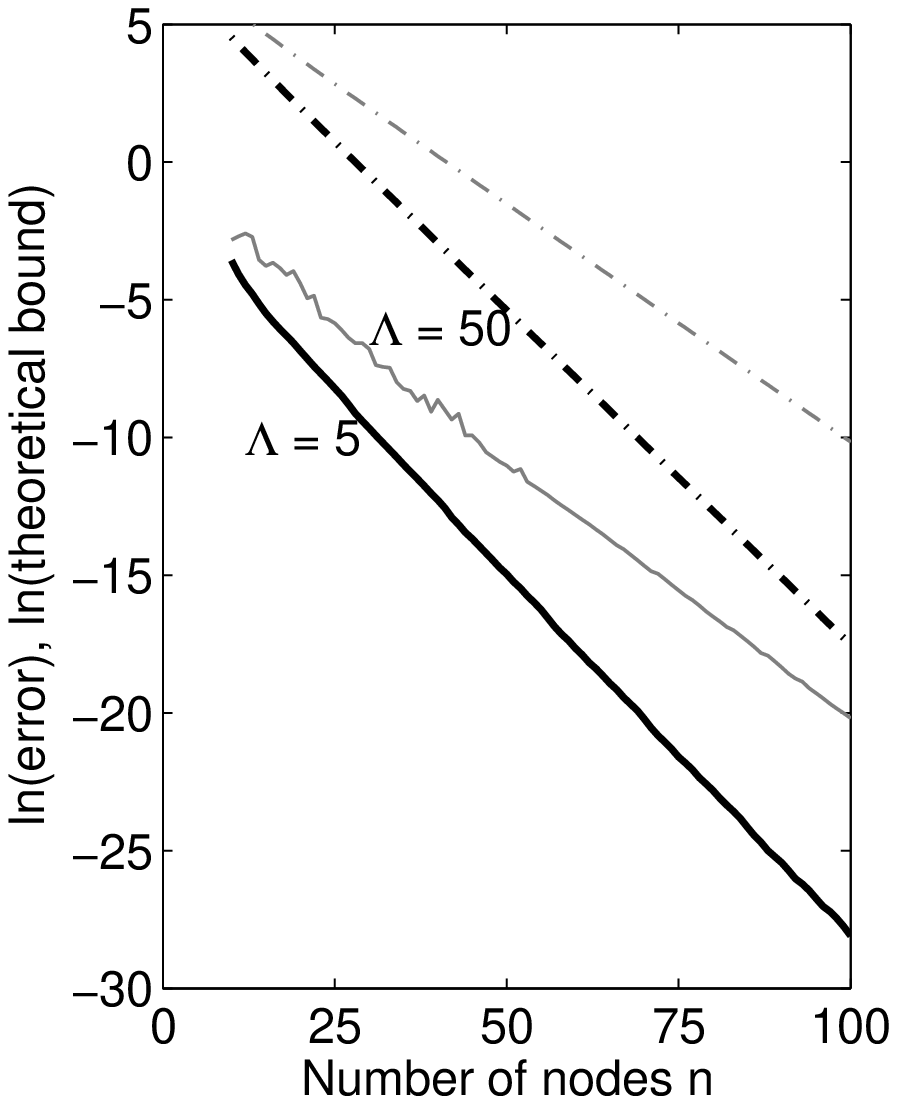, height=3in}
      \end{center}
      \caption{Natural logarithm of $\max_{t \in [t_0, \Lambda t_0]}\|u(t)-\bar{u}_n(t)\|$
      (continuous) and the theoretical estimate (dashed) versus $n$, for $u$ in
       Illustration~\ref{ilustracion2}. The gray lines correspond to $\Lambda = 50$
       and the black ones to $\Lambda=5$.}
      \label{figMLeffler}
    \end{figure}

  }
\end{ilustracion}

\begin{ilustracion}\label{ilustracion3}
{\rm We consider the inhomogenous heat equation on the unit square $\Omega =
  (0,1)^{2}$ with zero initial value and a convective heat flux at the boundary
  \begin{equation}\label{eq:heat2D}
    \left\{
      \begin{array}{lcl}
        u_t(t,x) &=& \Delta u(t,x) + f(x), \mbox{ for } x \in \Omega, \, t \ge 0, \\
        \partial_{\nu} u(t,x) &=& -u(t,x), \mbox{ for } x \in \partial\Omega, \, t
        \ge 0, \\ u(0,x) &=& 0, \mbox{ for } x \in \Omega,
      \end{array} \right.
  \end{equation}
where $f$ is the indicator function of the rectangle $R = [0.6,0.8]\times[0.2,0.8]$,
i.e. $f=1$ on $R$ and $f=0$ elsewhere.

Problem (\ref{eq:heat2D}) is semi-discretized in space by using linear finite elements
on a triangular grid. Denoting by $V_h \subset L^2(\Omega)$ the space of elements and
by $U_h(z)$ the Laplace transform of the semi-discrete solution $u_h(t)$, we get
$$ U_h(z) = \frac{1}{z} (z - \Delta_h)^{-1} P_hf,
$$ with $\Delta_h : V_h \to V_h$ the discrete Laplacian and $P_h$ the orthogonal
projection of $f$ onto $V_h$. Now, for fixed $h >0$, we try to approximate $u_h(t)$
by inverting $U_h(z)$. Notice that, since $\Delta_h$ is definite negative, certainly
$U_h(z)$ satisfies (\ref{desU}) for any $0 < \delta < \pi/2$ and $M =
1/(\mu_h\sin(\delta))$, with $-\mu_h$ the highest eigenvalue of $\Delta_h$. Notice
also that, working in coordinates relative to the standard basis of elements,
$U_h(z)$ is represented by a vector valued mapping ${\bf U}_h (z)$ satisfying
$$ zM_h {\bf U}_h (z) + S_h {\bf U}_h (z) = \frac{1}{z} {\bf f}_h ,
$$ where $M_h$ and $S_h$ stand for the mass and stiffness matrices and where ${\bf
f}_h$ is the vector formed by the scalar products of $f$ with the elements of the
basis. Thus, one evaluation of $U(z_k)$ at a given node $z_k$ requires the solution
of one linear system of the above form.

\begin{figure}[!ht]
  \begin{center}
    \begin{minipage}[c]{0.49\textwidth}
       \epsfig{file=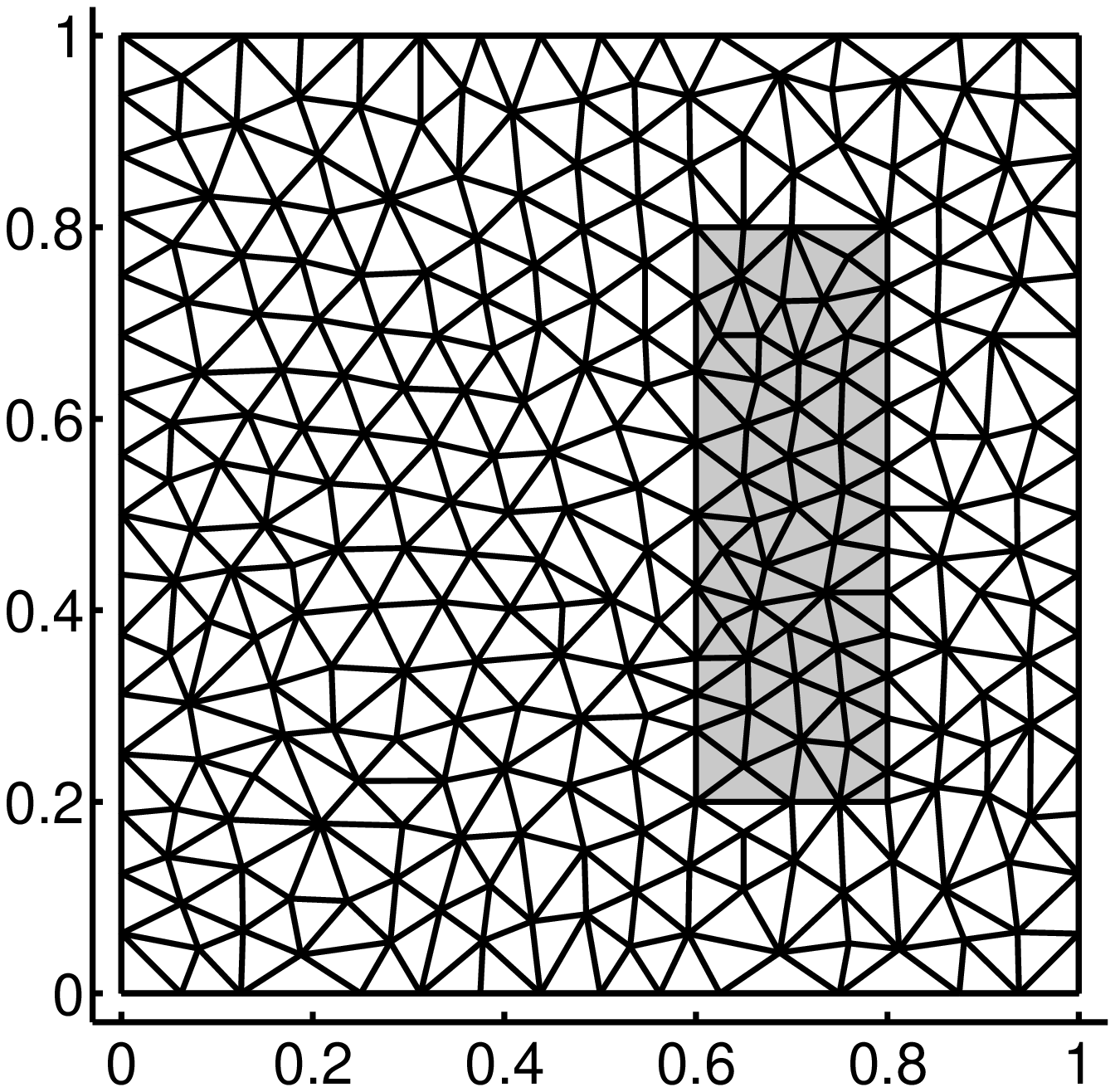, width=0.95\textwidth}
    \end{minipage}
    \begin{minipage}[c]{0.49\textwidth}
      \epsfig{file=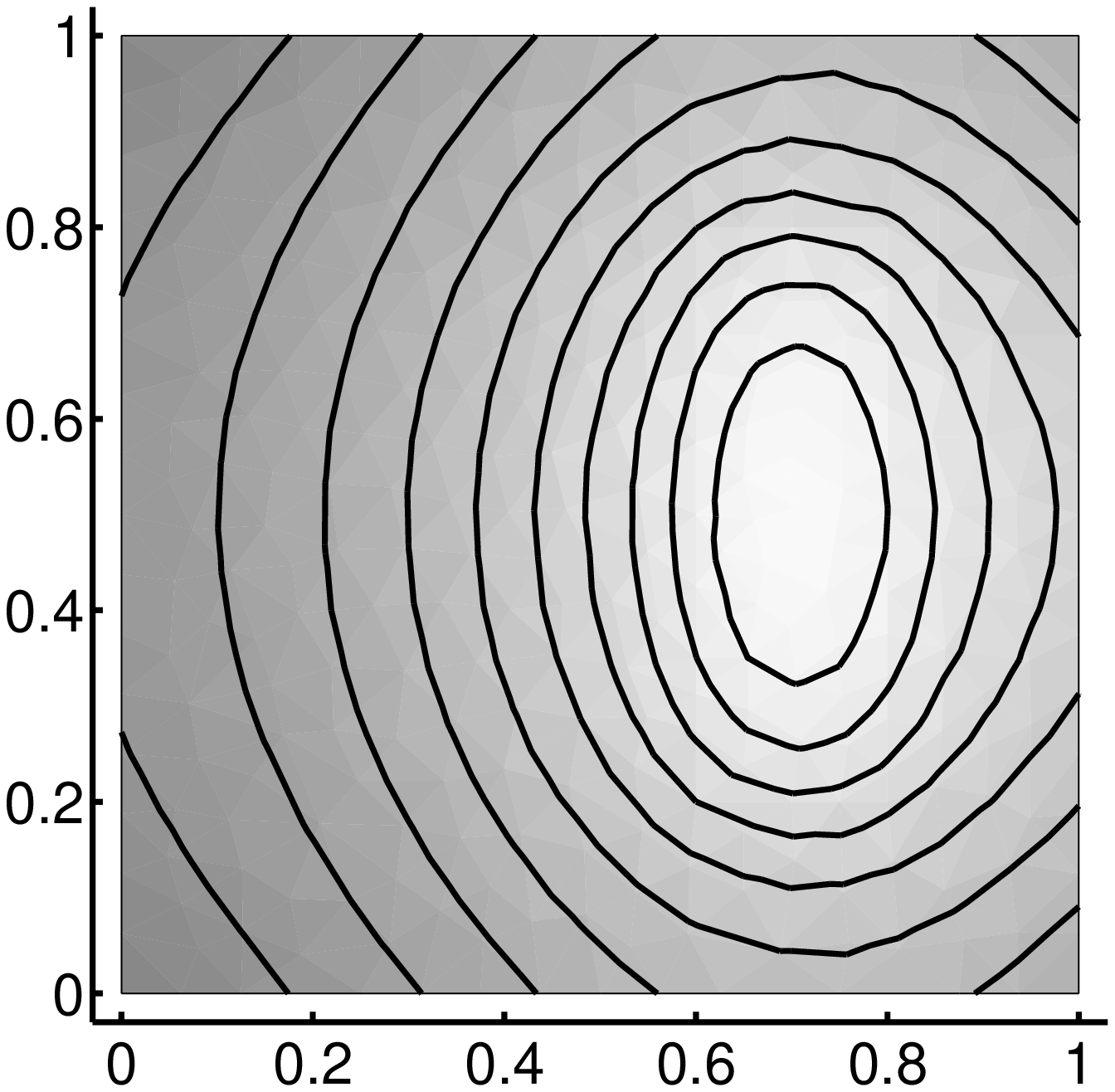, width=0.95\textwidth}
    \end{minipage}
  \end{center}
  \caption{Left: Mesh of $\Omega$, with the set $R$ indicated by dark-gray.
    Right: Temperature distribution at $t=0.5$ in false-color representation.  (white
    corresponds to temperature $1$ and black to $0$)}
  \label{figheat_geo_cont}
\end{figure}

In the experiment we generate a mesh, shown in the left of
Fig~\ref{figheat_geo_cont}, with 542 triangles by means of the mesh generator
Triangle~\cite{Shewchuk}. Linear systems are solved using MATLABs sparse LU
factorization UMFPACK. Since $u_h(t)$ is unknown, the errors are estimated in the
$L^2(\Omega)$-norm with respect to a reference solution $\bar u_{h,500}(t)$ obtained
with 500 nodes. In the absence of precise information about $\rho$, both for this
reference solution and for the rest of the approximations $\bar{u}_{h,n}(t)$ to
$u_h(t)$, we tune $\theta = 1-1/n$, as indicated in Section~\ref{secparametros}. In
Fig.~\ref{figheat}, for the parameters $\alpha = 0.7$, $d = 0.6$, $t_0=0.01$ and
$\theta=1-1/n$, we plot $\ln \max_{t \in [t_0,\Lambda t_0]}
\|\bar{u}_{h,500}(t)-\bar{u}_{h,n}(t)\|$ against $n$, for $\Lambda = 5,\, 50$. This
plot shows the predicted behavior.

\begin{figure}[!ht]
  \begin{center}
      \epsfig{file=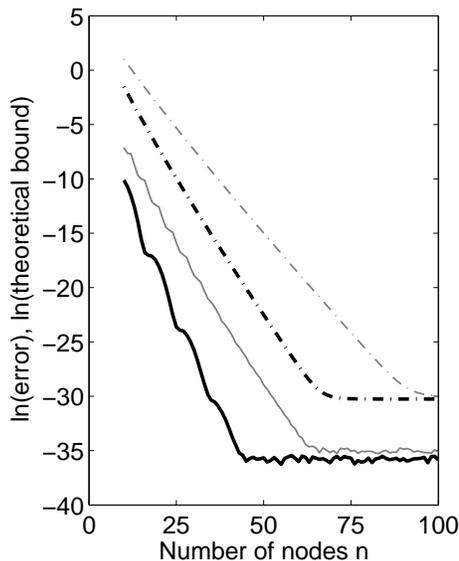, height=3in}
  \end{center}
  \caption{Left: Natural logarithm of $\max_{t\in [t_{0}, \Lambda
      t_{0}]}\|u(t)-\bar{u}_n(t)\|$ (continuous) and the theoretical estimate
    (dashed) versus $n$, for $u$ in Illustration~\ref{ilustracion3}.  The gray lines
    correspond to $\Lambda = 5$ and the black ones to $\Lambda = 50$.}
  \label{figheat}
\end{figure}

  }
 \end{ilustracion}

\begin{ilustracion}\label{ilustracion4}
{\rm We consider again the Laplace transform $U(z) = 1/(1+z)$ of the exponential
function $u(t) = e^{-t}$ as in Illustration~\ref{ilustracion1}. The values of
$\alpha,\, d$ and $t_0$ are again $0.7,\ 0.6$ and $1$, respectively.

We add on purpose perturbations of maximum size $10^{-4}$ to the evaluations of $U$
at the required nodes. Thus, we use (\ref{invnumericareal}) with
$$ U_k = U(z_k) + \eta_k, \qquad -n \le k \le n,
$$ with $|\eta_k| \le \rho = 10^{-4}$. Now we try to approximate $u(t)$ without using
the available information about $\rho$. In this situation, as explained in
Section~\ref{secparametros}, we take $\theta=1-1/n$.

In fact, we compare two types of perturbations:

We first generate complex, random, independent perturbations $\eta_k$, in such a way
that $|\eta_k|$ and $\arg (\eta_k)$ are uniformly distributed on $[0,10^{-4}]$ and
$[0,2\pi]$, respectively. In Fig.~\ref{figpertubacion} (left), we show the resulting
actual error, which behaves much better than predicted by (\ref{cotanooptimal}). The
explanation is that cancellations are likely compensating the effects of the
independent random perturbations. A finer analysis of the observed behavior is out of
the scope of the present paper.

Secondly, for each $-n\leq k\leq n$, we consider the perturbation
$$ \eta_k = 10^{-4} \exp(-i \arg (\omega_k(t_0))),$$ with $\omega_k(t_0)$ defined in
(\ref{invnumericareal}). These perturbations correspond to the worst possible case in
(\ref{abserrorevaluacion}), for $t=t_0=1$. Now, the resulting actual error, plotted
in Fig.~\ref{figpertubacion} (right), fits quite well with (\ref{cotanooptimal}).

\begin{figure}
      \begin{center}
        \epsfig{file=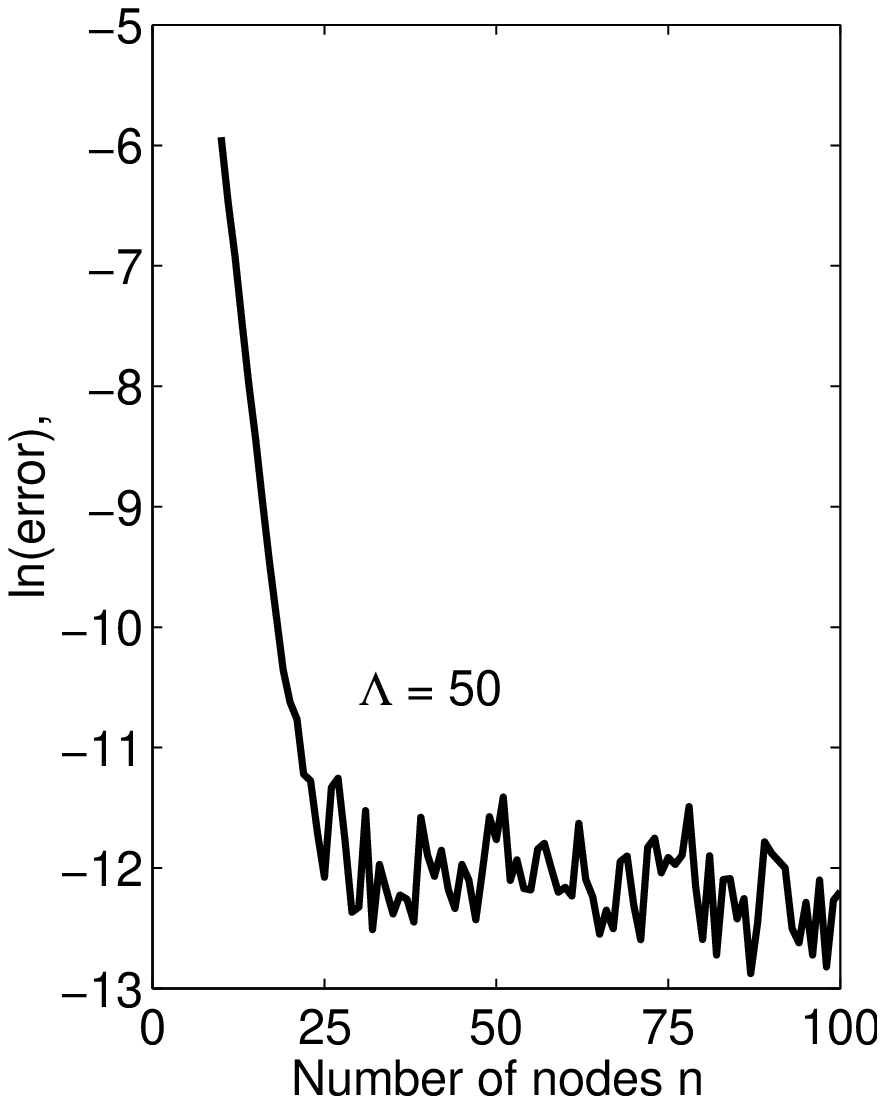, height=3in}
        \epsfig{file=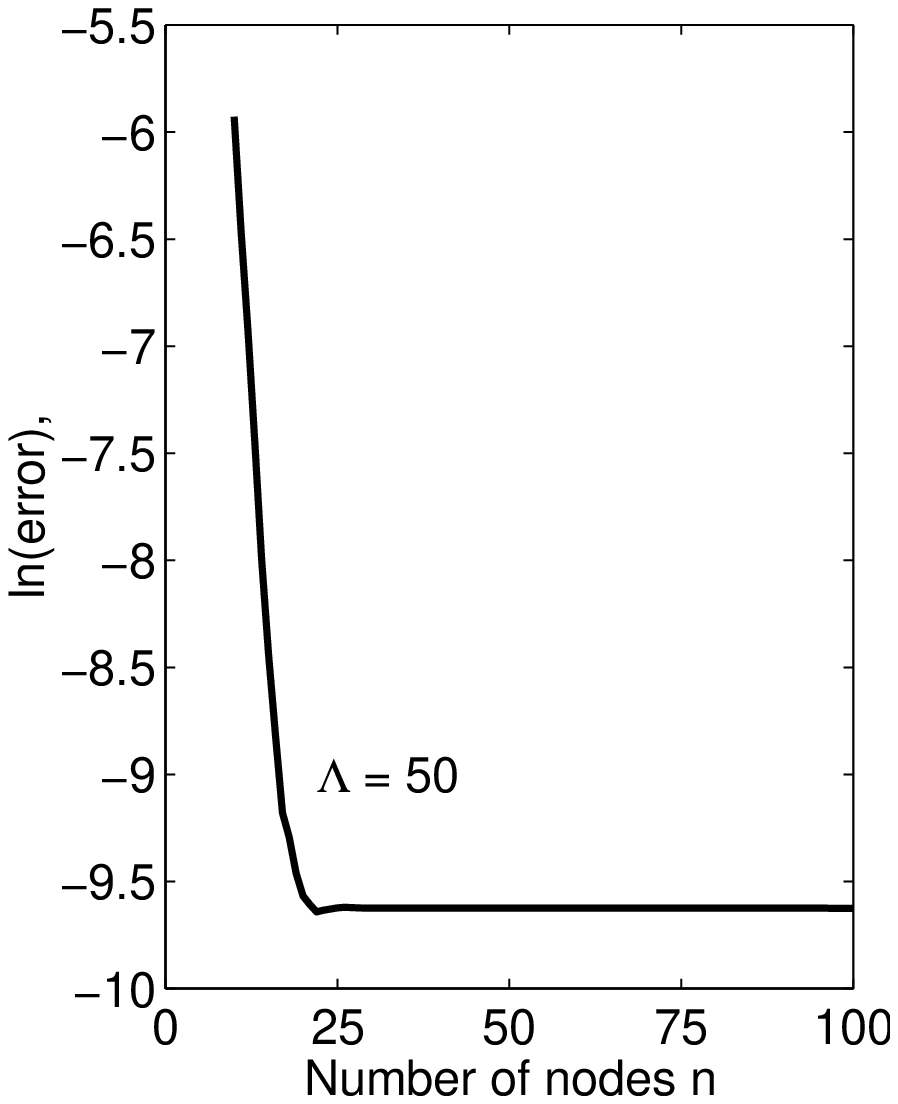, height=3in}
      \end{center}
      \caption{$\ln \max_{t\in [t_0, \Lambda t_0]} (\|u(t)-\bar{u}_n(t)\|)$
      versus $n$, for $u$ in Illustration~\ref{ilustracion4} with $\theta=1-1/n,\
      \alpha=0.7,\ d=0.6$, $t_0=1$ and $\Lambda = 50$. Left: Random perturbation.
      Right: Worst case perturbation.}
      \label{figpertubacion}
    \end{figure}

  }
\end{ilustracion}

\bibliographystyle{amsplain}

\end{document}